\documentclass[12pt,notitlepage,twoside]{article}
\usepackage{latexsym} 
\usepackage{amsmath} 
\usepackage{amsfonts} 
\usepackage{amssymb} 
\renewcommand{\a }{\alpha } 
 
\renewcommand{\d}{\delta }

\newcommand{\D }{\Delta }
 
\newcommand{\e }{\varepsilon }

\renewcommand{\l }{\lambda } 
\renewcommand{\L }{\Lambda } 
 
\newcommand{\n }{\nabla } 
\newcommand{\var }{\varphi }

\newcommand{\Sig }{\Sigma} 
\renewcommand{\t }{\tau } 
 
\renewcommand{\o }{\omega }

\newcommand{\ov}{\overline}

\newcommand{\be}{\begin{equation}} 
\newcommand{\ee}{\end{equation}} 
\newenvironment{pf}{\noindent{\bf Proof.}\enspace}{%\rule{2mm}{2mm}
\hfill$\Box$\medskip} 
\newenvironment{pfn}[1]{\noindent{\bf Proof of {#1}\enspace}}{%\rule{2mm}{2mm}
\hfill$\Box$\medskip} 
\newcommand{\R}{\mathbb{R}} 
\newcommand{\N}{\mathbb{N}}

\newtheorem{thm}{Theorem}[section] 
\newtheorem{pro}[thm]{Proposition}
\newtheorem{lem}[thm]{Lemma}

\newtheorem{cor}[thm]{Corollary}
 
\numberwithin{equation}{section}

\author{{ {\bf M. Ben Ayed },  {\bf K.  El Mehdi} 
 \& {\bf M.  Ould Ahmedou}}}

 \title { \Large \textbf{The Scalar Curvature Problem on the Four \\
  Dimensional  Half Sphere}}

\begin{document}

\date{ }

\maketitle

{\footnotesize

\noindent 
{\bf Abstract.}
In this paper, we consider the problem of prescribing the scalar curvature under minimal boundary conditions on the standard four dimensional half sphere. We provide an Euler-Hopf type criterion for a given function to be a scalar curvature to a metric conformal to the standard one. Our proof involves the study of critical points at infinity of the associated variational problem.
}

\section{Introduction and the Main Results }

This paper is devoted to some nonlinear problem arising from conformal geometry.
Precisely, Let  $(M^n,g)$ be a $n$-dimensional  Riemannian manifold with boundary,
$n\geq 3$, and let $\tilde{g}=u^{4/(n-2)}g$, be a conformal metric to $g$, 
where $u$ is a smooth positive function, then  
the scalar curvatures $R_g$, $R_{\tilde{g}}$ and the mean curvatures  $h_g$, $h_{\tilde{g}}$, with respect to $g$ and $\tilde{g}$ respectively, are related by the following equations. 
\begin{eqnarray*}
(P_1)\quad \left\{
\begin{array}{ccccc} 
-c_n\D _gu+R_gu&=&R_{\tilde{g}}u^{\frac{n+2}{n-2}}&\mbox{ in }& M\\
\frac{2}{n-2}\frac{\partial u}{\partial \nu}+h_g u &=&
h_{\tilde{g}}u^\frac{n}{n-2}&
\mbox{ on }& \partial M
\end{array}
\right.
\end{eqnarray*}
where $c_n=4(n-1)/(n-2)$ and $\nu$ denotes the outward normal vector with 
respect to the metric $g$.\\
In view of $(P_1)$,  the following problem naturally arises   : given two functions  $K: M\to \R$ and $H: \partial M\to \R$, does exist a metric $\tilde{g}$
conformally equivalent to $g$ such that $R_{\tilde{g}}=K$ and 
$h_{\tilde{g}}=H$? \\
From  equations $ (P_1)$, the problem is equivalent to finding a smooth 
positive solution $u$ of the following equation
\begin{eqnarray*}
(P_2)\quad \left\{
\begin{array}{ccccc} 
-c_n\D _gu+R_gu&=& K u^{\frac{n+2}{n-2}}&\mbox{ in }& M \\
\frac{2}{n-2}\frac{\partial u}{\partial \nu}+h_g u &= & 
 H u^\frac{n}{n-2}&
\mbox{ on }& \partial M
\end{array}
\right.
\end{eqnarray*}
This problem has been studied in earlier works   (see \cite{ALM} \cite{BEO},\cite{C},\cite{DMO}  \cite{E1}, \cite{E2}, \cite{Escobar}, \cite{HL1}, \cite{HL2}, \cite{Hebey}, \cite{Li} and the references therein).\\
In this paper we consider the case of the standard four dimensional half sphere under minimal boundary conditions. We are thus reduced to look for positive solutions of the following problem
\begin{eqnarray*}
(1)\quad \left\{
\begin{array} {ccccc}
L_gu := -\D _gu+ \frac{n(n-2)}{4}u & = & K u^{\frac{n+2}{n-2}} &\mbox{ in }& S_+^n\\
\frac{\partial u}{\partial \nu}&=&0&  
\mbox{ on }& \partial S_+^n
\end{array}
\right.
\end{eqnarray*}
where  $n=4$ and where $g$ is the standard metric of $S_+^n =\{x\in \R^{n+1}\, /\, |x|=1,\, x_{n+1}>0\}$.\\

This problem has been studied by
 Yanyan Li \cite{Li}, and Djadli-Malchiodi-Ould Ahmedou  \cite{DMO}, on the three dimensional standard half
sphere. Their method involves a fine blow up analysis of some subcritical 
approximations and the use of the topological degree tools. In a previous work \cite{BEO}, we gave some topological conditions on $K$ to prescribe the scalar curvature under minimal boundary conditions on half spheres of dimension bigger than or equal $4$ . In this paper, which is a continuation of \cite{BEO}, we single out the four dimensional case to give more existence results, in particular a Hopf type formula reminiscent to Bahri-Coron \cite{BaC} formula for the scalar curvature problem on $S^3$, see also \cite{BCCH} \cite{cgy}.\\
Notice that, Problem $(1)$ is in a natural way  related to the well-known Scalar Curvature Problem on $S^n$
$$
(2) \qquad -\D_g u + \frac{n(n-2)}{4}u = K u^{(n+2)/(n-2)} \qquad \mbox{in } S^n$$
to which much works have been devoted (see \cite{agp},\cite{aujfa},\cite{ab},\cite{B1},  \cite{B2}, \cite{BCCH},\cite{cl},\cite{cgy},\cite{cy1},\cite{C},\cite{es}, \cite{hv}, \cite{y1}, \cite{[SZ]} and the references therein ). As for Problem (2), also for
problem (1) there are topological obstructions for
existence of solutions, based on Kazdan-Warner type conditions,
see \cite{bp} . Hence
it is not expectable to solve problem (1) for all the
functions $K$ , and it is natural to impose some conditions
on them.\\
Regarding Problem $(2)$, Bahri \cite{B2} observed that a new phenomenon appears in dimension $n\geq 5$ due to the fact that the self interaction of the functions failing the Palais-Smale condition dominates the interaction of two of those functions. In the dimensional  three, the reverse happens(see \cite{BaC}). In dimension $4$, we have a balance phenomenon, that is, the self interaction and the interaction are of the same size (see \cite{BCCH}).\\
For Problem $(1)$, Djadli-Malchiodi-Ould Ahmedou \cite{DMO} showed that such a  balance phenomenon appears in $S^3_+$. Such a result suggests that there is a dimension break between Problem $(1)$ and Problem $(2)$. In this paper, we prove that a dimension's  break is  not always true. Precisely, we show that $S^4_+$ behaves like $S^4$ under some conditions on the behavior of the function $K$ on the boundary $\partial S^4_+ $. \\
%Our approch borrows some of the  ideas developped in \cite{BCCH} where the pro%blem of prescribing the scalar curvature on four closed manifolds was studied %using  some  Euler-Poincar\'e Characteristic argument.\\
In order to state our results, we need to introduce some notations and the assumptions that we are using in our results.\\
Let $G$ be the Green's function of $L_g$ on $S^4_+$ and $H$ its regular part defined by 
\begin{eqnarray*}
\begin{cases}
G(x,y)=(1-\cos(d(x,y)))^{-1}+ H(x,y),\\
\D H =0 \mbox{ in } S^4_+,\quad \partial G/\partial\nu =0 \mbox{ on } \partial S^4_+
\end{cases}
\end{eqnarray*} 
 Let $K$ be a $C^3$ positive Morse function on $\ov{S^4_+}$. We say that  that $K$ satisfies the condition $(\textbf{C})$ if 
\begin{eqnarray*}
\begin{cases}
\mbox{If } a \in S^4_+ , \n K(a)= 0, \qquad \mbox{then } \frac{-\D K(a)}{3K(a)}-4 H(a,a) \neq 0, \\

\mbox{If } a \in \partial S^4_+ , \n_T K(a)= 0, \qquad \mbox{then }\frac{\partial K}{\partial \nu}(a) < 0,
\end{cases}
\end{eqnarray*}
where $ \n_T K$ denotes the tangitial part of the gradient of $K$

% only nondegenerate critical points $y_1,...,y_m$ such that
%$$\frac{-\D K(y_i)}{3K(y_i)} +4 H(y_i,y_i) \neq 0, \qquad \mbox{ for } \, i=1,%...,N
%$$
For  sake of simplicity, we assume that $y_1,...,y_l$ ($l\leq N$) are all the critical points of $K$ for which $(-\D K(y_i)/(3K(y_i))) - 4 H(y_i,y_i) >0$, for $ i=1,...,l$. 
%We also assume that all the critical points of  
%$K_1=K_{/\partial S^4_+}$ are $z_1,...,z_{m'}$, and satisfy 
%$$
%\frac{\partial K}{\partial \nu}(z_i) < 0, \qquad \mbox {for } i=1,...,m'
%$$
For $s\in \{1,...,l\}$ and for any $s$-tuple $\tau_s=(i_1,...,i_s) \in \{1,...,l\}^s$ such that $i_p\neq i_q$ if $p\neq q$, we define a Matrix $M(\tau_s)=(M_{pq})_{1\leq p,q\leq s}$, by
$$
M_{pp}=\frac{-\D K(y_{i_p})}{3K(y_{i_p})^2} - 4\frac{ H(y_{i_p},y_{i_p})}{K(y_{i_p})}, \qquad M_{pq} =-\frac{4G(y_{i_p},y_{i_q})}{\left(K(y_{i_p})K(y_{i_q})\right)^{1/2}} \quad \mbox{for } p\neq q,
$$
and we denote by $\rho (\tau _s)$ the least eigenvalue of $M(\tau _s)$.\\
We then have the following results

\begin{thm}\label{t:11}
Assume that $K$ satisfies condition $(\textbf{C})$, and assume further that for any $s\in \{1,...,l\}$,
$M(\tau_s)$ is nondegenerate.\\
 If
$$
1\ne \sum_{s=1}^l\sum_{\tau_s=(i_1,...,i_s)/\rho(\tau_s)>0}(-1)^{s-1-
\sum_{j=1}^sk_{i_j}},
$$
where $k_{i_j}= \mbox{index }\left(K, y_{i_j}\right)$, denotes the Morse Index of $K \mbox{at }  y_{i_j}$ , then Problem $(1)$ has a solution.
\end{thm}
\begin{thm}\label{t:12}
Assume that $K$ satisfies condition $(\textbf{C})$, and assume further that:

$$
\mbox{If } a \in S^4_+ , \n K(a)= 0,  \mbox{we have } -\frac{\D K(a)}{3K(a)} \leq \frac{1}{d^2(a, \partial S^4_+)}  ,
$$
 then Problem $(1)$ has a solution.
\end{thm}

The proof of the above results involves the construction of a special pseudogradient for the associated variational problem for which the Palais Smale condition is satisfied along the decreasing flow lines , as long as these flow lines do not enter the neighbourhood of a finite numbers of critical points of $K$ such that  the related matrix $M(\t)$ is positive definite.
Moreover along the flow lines of such a pseudogradient there can be only finitely many blow up points. Furthermore if some blow up points are close and the interactions between them is large, then the flow lines starting from there will enter the zone with at least one less blow up points. Using  such a pseudogradient, a Morse Lemma at infinity is performed and an  Euler Poincar\'e characteristic argument allows us to derive the existence of solution for this problem.\\
Our proof goes along the method initiated by Bahri-Coron \cite{BaC}, see also \cite{BCCH}, however in our case the presence of the boundary makes the analysis more involved. Moreover it turns out that the interaction of the bubbles and the boundary creates a  phenomenon of new type which is not present in the sphere's case. \\

We organize the remainder of the present paper as follows. In section 2, we set up the variational structure and  give some careful expansion of the Euler functional $J$ associated to $(1)$ . In section 3, we perform the construction of a pseudogradient of $J$ whose zeros are the critical points at infinity of $J$. The last section is devoted to the proofs of our results.\\

\noindent {\bf Acknowledgment.} The authors would like to thank Professor Abbas Bahri for his encouragememt and constant support over the years. The third Author (M.O.A) is grateful to Professor Reiner Sch\"{a}tzle for his support and interest in this work.

\section { General Framework and Expansion of the functional}
\mbox{}
Problem (1) has a variational structure, the functional being 
$$
J(u)= \frac {\int_{S_+^4}|\n u|^2+2\int_{S_+^4}u^2}
{\left(\int_{S_+^4}Ku^4\right)^{\frac{1}{2}}},
$$
defined on the unite sphere of $H^1(S^4_+)$ equipped with the norm 
$$
||u||^2=\int_{S_+^4}|\n u|^2+ 2 \int_{S_+^4}u^2.
$$
Problem (1) is equivalent to finding the critical points of $J$ to the constraint $u\in \Sig^+$, where 
$$\Sigma ^+=\{u \in \Sigma \, / \, u\geq 0\}, \quad \Sig=\{u\in H^1(S^4_+)/\, \, ||u||=1\}
$$
The Palais-Smale condition fails to be satisfied for $J$ on $\Sigma ^+$. To describe the sequences failing the Palais-Smale condition, we 
need to introduce some notations.\\
For $a\in \overline{S_+^4}$ and $\l >0$, let 
$$
\d _{a,\l }(x)=\frac {\l}{\l ^2+1+(1-\l ^2)
\cos d(a,x)},
$$
where $d$ is the geodesic distance on $(\ov{S^4_+},g)$. This function satisfies the following equation  
$$
-\D \d _{a,\l } + 2\d _{a,\l }=
8\d _{a,\l }^3, \quad \mbox{ in } S_+^4
$$
Let $\var_{(a,\l)}$ be the function defined on $ S_+^4$ and satisfying 
$$
-\D \var _{(a,\l) } + 2\var _{(a,\l) }= -\D \d _{a,\l } + 2\d _{a,\l }  \mbox{ in } S_+^4, \quad \frac{\partial\var_{(a,\l)}}{\partial\nu}=0  \mbox{ on } \partial S_+^4
$$
Regarding $ \var_{(a,\l)}$ we prove the following lemma
\begin{lem}\label{l:21}
For $a\in \partial S_+^4$, we have $(\partial\d_{(a,\l)})/(\partial\nu)=0$ and therefore $\var_{(a,\l)}=\d_{(a,\l)}$. For $a\notin \partial S_+^4$, we have 
$$
\var_{(a,\l)}=\d_{(a,\l)}+\frac{H(a,.)}{\l}+f_{(a,\l)},
$$  where $ f_{(a,\l)}$ satisfies  
$$|f_{(a,\l)}|_{L^\infty}\leq \frac{c}{\l^3d^4}, \quad \l\frac{\partial f}{\partial\l}=O\left( \frac{1}{\l^3d^4}\right)$$
 and where $d=d(a,\partial S^4_+)$.
\end{lem}
\begin{pf}
Using a stereographic projection, we are led to prove the corresponding estimates  on $\R^4_+$. We still   denote by $G$ and $H$ the Green's function and its regular part of Laplacian on $\R^4_+$ under Neumann boundary conditions. In this case, we have 
$$\d_{(a,\l)}=\frac{\l}{1+\l^2|a-x|^2}\quad  \mbox{ and } \quad H(a,x)=\frac{1}{|\bar{a}-x|^2},$$
where $\bar{a}$ is the symmetric of $a$ with respect to $\partial \R^4_+$.\\
Observe that, for $\theta=\d_{(a,\l)}+\var_{(a,\l)}-H(a,.)/\l$, we have 
$$\D \theta =0 \mbox{ in } \R^4_+ , \quad \frac{\partial \theta}{\partial \nu}=\frac{\partial \d}{\partial \nu}+\frac{1}{\l}\frac{\partial H}{\partial \nu}= O\left(\frac{1}{\l^3d^5}\right)$$
Thus, using the Green's formula, we derive
$$ \theta(y)=c\int_{\partial\R^4_+}G \biggl(\frac{\partial \d}{\partial \nu}+\frac{1}{\l}\frac{\partial H}{\partial \nu}\biggr)\leq \frac{c'}{\l^3d_a^2} \int_{\partial\R^4_+}G\frac{1}{|a-x|^3 }, $$ 
where $d_a$ is the distance of $a$ to the boundary. But $G$ satisfies 
$$\int_{\partial \R^4_+}G(x,y)\frac{1}{|a-x|^3}=O\left(\frac{1}{d_a ^2}\right).$$
 Thus, the first estimate follows. The second estimate can be proved by the same way. 
\end{pf}\\
Now, for $\e >0$ and $p\in \N^*$, let us define  
\begin{align*}
V(p,\e )=& \{u\in \Sig /\exists a_1,...,a_p \in \overline{S_+^4}, \exists 
\l _1,...,\l _p >0, \exists \a _1,...,\a _p>0 
 \mbox{ s.t. } ||u-\sum_{i=1}^p\a _i\d _i||<\e, \\
 &  |\frac 
{\a _i ^2K(a_i)}{\a _j ^2K(a_j)}-1|<\e ,
 \l _i>\e ^{-1}, \e _{ij}<\e \mbox{ and }\l _id_i<\e \mbox{ or }
\l _id_i >\e ^{-1}\},
\end{align*}
where $\d _i=\d _{a_i,\l _i}$, $d_i=d(a_i,\partial S_+^4)$ and 
$\e _{ij}^{-1}=\l _i/\l _j +\l _j/\l _i + \l _i\l _j(1-\cos d(a_i,a_j))/2$.\\
The failure of Palais-Smale condition can be described, following the ideas introduced in \cite{BC}, \cite{L},  \cite{S} as follows:

\begin{pro}\label{p:22} 
Assume that $J$ has no critical point in $\Sigma ^+$ and let $(u_k)\in 
\Sigma ^+$ be a sequence such that $J(u_k)$ is bounded and $\n J(u_k)\to 0$.
Then, there exist an integer $p\in \N^*$, a sequence $\e _k>0$ ($\e _k\to 0$)
and an extracted subsequence of $u_k$, again denoted $(u_k)$, such that 
$u_k\in V(p,\e _k )$.
\end{pro}

If a function $u$ belongs to $V(p,\e)$, we assume , for the sake of simplicity,that  $\l_id_i < \e$ for $i\leq q$ and $\l_id_i > \e^{-1}$ for $i > q$. We consider the following minimization problem for $u\in V(p,\e)$ with $\e$ small 
\begin{eqnarray}\label{e:51}
\min\{||u-\sum_{i=1}^q\a _i\d _{(a_i,\l_i)}-\sum_{i=q+1}^p\a _i\var_{(b_i,\l_i)} ||,\, \a _i>0,\, \l _i>0,\, a_i\in 
\partial S_+^4 \mbox{ and } b_i\in S_+^4\}
\end{eqnarray}
We then have the following proposition which defines a parametrization of the set $V(p,\e )$. It follows from corresponding  statements in \cite{B2}, \cite{BC'}, \cite{R}.

\begin{pro}\label{p:23} 

For any $p\in \N^*$, there is $\e _p>0$ such that if $\e <\e _p$ and $u\in V(p,\e )$, the minimization problem \eqref{e:51}  
has a unique solution (up to permutation). In particular, we can write 
$u\in V(p,\e )$ as follows 
$$
u=\sum_{i=1}^q\bar{\a }_i\d _{(\bar{a}_i,\bar{\l }_i)}+\sum_{i=q+1}^p\bar{\a }_i\var _{(\bar{a}_i,\bar{\l }_i)}+ v,
$$
where $(\bar{\a }_1,...,\bar{\a }_p,\bar{a}_1,...,\bar{a}_p,\bar{\l }_1,...,
\bar{\l }_p)$ is the solution of \eqref{e:51} and $v\in H^1(S_+^n)$ such that 
$$(V_0)\qquad ||v|| \leq \e, \quad
(v,\psi)=0 \mbox{ for } \psi\in \bigg\{\d_i,\frac{\partial\d_i}{\partial\l_i},\frac{\partial\d_i}{\partial a_i},\var_j,\frac{\partial\var_j}{\partial\l_j},\frac{\partial\var_j}{\partial a_j}/ \, \, i\leq q , \, j>q\bigg\} 
$$
\end{pro}
We also have the following proposition whose proof is similar, up to minor modification to corresponding statements in \cite{B1} (see also \cite{R}) 

\begin{pro}\label{p:24}
There exists a $C^1$ map which, to each $(\a_1,...,\a_ p,a_1,...,a_ p,\l_1,..., \l_ p)$ such that 
$ \sum_{i=1}^p\a _i\d _i \in V(p,\e )$ with small $\e $, 
associates $\ov{v}=\ov{v}_{(\a_i,a_i,\l_i )}$ satisfying 
$$
J\left(\sum_{i=1}^q\a _i\d _i +\sum_{i=q+1}^p\a_i\var_i+\ov{v}\right)= \min\bigg\{
J\left( \sum_{i=1}^q\a _i\d _i+\sum_{i=q+1}^p\a_i\var_i +v\right) , \, v \mbox{ satisfies } (V_0)\bigg\}.
$$
Moreover, there exists $c>0 $ such that the following holds
$$
||\ov{v}||\leq c \left(\sum_{i\leq q}\frac{1}{\l _i}+\sum_{i>q}\frac{|\n K(a_i)|}{\l_i}+\sum_{i>q}\frac{1}{(\l_id_i)^2}+\sum_{k\ne r}\e _{kr} 
(log (\e _{kr}^{-1}))^{1/2}\right).
$$
\end{pro}
Next, we are going to give an useful expansion of functional $J$ and its gradient in $V(p,\e)$. 

\begin{pro}\label{p:J}
For $\e>0$ small enough and $u=\sum_{i=1}^p\a_i\var_{(a_i,\l_i)}\in V(p,\e)$, we have the following expansion
\begin{align*}
J(u)= & \frac{8S^{1/2}\sum\a_i ^2}{\left(\sum\a_i ^4K(a_i)\right)^{1/2}}\left(1+\frac{\o_3}{8S}\left(\sum K(a_i)^{-1}\right)^{-1}\left(\sum\biggl(\frac{-\D K(a_i)}{3\l_i ^2K(a_i)^2}-\frac{4H(a_i,a_i)}{\l_i ^2K(a_i)}\biggr)\right.\right.\\
 &\left.\left. -\sum_{i\ne j}\frac{2}{\left(K(a_i)K(a_j)\right)^{1/2}}\left(\e_{ij}+\frac{2H(a_i,a_j)}{\l_i\l_j}\right)\right)+o\left(\sum \e_{kr}+\frac{1}{(\l_kd_k)^2}\right)\right),
\end{align*} 
where $$ S= \int_{\R^4} \frac{dx}{(1 +|x|^2 )^4}$$
\end{pro}
\begin{pf}
We need to estimate 
$$
N=||u||^2 \mbox{ and } D^2=\int_{S^4_+}K u^4
$$
Observe that 
\begin{eqnarray}\label{e:1'}
N= \sum\a_i ^2||\var_i||^2+\sum_{i\ne j}\a_i\,\a_j\,(\var_i,\var_j)
\end{eqnarray}
As in \cite{BEO}, we  have 
\begin{eqnarray}\label{e:2'}
||\var_i||^2=8\, S + 8 \omega \frac{H(a_i,a_i)}{2\l_i ^2}+o\left(\frac{1}{(\l_id_i)^2}\right)
\end{eqnarray}
We also have, for $i\ne j$
 \begin{eqnarray}\label{e:3'}
(\var_i,\var_j)=2\o_3\left(\e_{ij}+\frac{2H(a_i,a_j)}{\l_i\l_j}\right)+o\left(\sum \frac{1}{(\l_kd_k)^2}\right)
\end{eqnarray}
For the denominator, we have
\begin{eqnarray}\label{e:4'}
D^2 =\int_{S_+^4}K \left(\sum\a_i\var_i\right)^4=\sum\a_i ^4\int_{S_+^4} K \var_i ^4+4\sum_{i\ne j}\a_i^3\a_j\int_{S_+^4}K\var_i ^3\var_j+o(\sum\e_{kr})
\end{eqnarray}
Observe that 
\begin{eqnarray}\label{e:5'}
\int_{S_+^4}K\var_i ^4=K(a_i)S+\frac{\o_3\D K(a_i)}{12\l_i}+2\o_3K(a_i)\frac{H(a_i,a_i)}{\l_i ^2}+o(\frac{1}{(\l_id_i)^2})
\end{eqnarray}
We also have, for $i\ne j$
 \begin{eqnarray}\label{e:6'}
\int_{S_+^4}K\var_i ^3\var_j=\frac{\o_3}{4}K(a_i)\left(\e_{ij}+\frac{2H(a_i,a_j)}{\l_i\l_j}\right)+o\left(\sum\frac{1}{(\l_kd_k)^2}\right)
\end{eqnarray}
Using \eqref{e:1'},...,\eqref{e:6'}, the result follows.
\end{pf}
\begin{pro}\label{p:25} (see \cite{BEO})
For $ u=\sum_{i\leq q}\a _i\d _i+\sum_{j=q+1}^p\a _j\var _j \in V(p,\e )$, we 
have the following expansions 
\begin{align*}
(\n J(u),\l _i\partial\d _i/\partial\l _i)= & c_1J(u)\sum_{i\leq q, j\ne i}\a_k
\l _i\frac{\partial \e _{ik}}{\partial \l _i}+c_2J(u)^3\frac{\a_i ^3}{\l _i}\frac{\partial K}{\partial \nu}(a_i)\\
 & +O\biggl(\sum_{k\leq q}\frac{1}{\l _k^2}+\sum_{r>q,k\leq q}\e_{kr}\biggr)+ 
o( \sum_{k,r\leq q}\e _{kr})\\
(\n J(u),\l _i^{-1}\partial\d _i/\partial a _i)= & c_3\a_iJ(u)e_4\biggl(c_4(1-J(u)^2\a_i ^2K(a_i))+J(u)^2\a_i ^2\frac{c_5}{\l_i}\frac{\partial K(a_i)}{\partial \nu}\biggr)\\
 & -c_5J(u)\sum_{k\leq q}\a_i\frac{1}{\l_i}\frac{\partial\e_{ik}}{\partial a_i}(1+o(1))-4J(u)^3\a_i ^3\frac{c_6}{\l_i}\n_TK(a_i)\\
 & +o(\sum_{k,r\leq q}\e_{kr})+O\biggl(\sum_{k\leq q}\frac{1}{\l_k^2}+\sum_{k\leq q,r>q}\e_{kr}\biggr)\\
(\n J(u),\d _i)= c_7J(u) \a _i & (1-J(u)^2\a _i^2K(a_i)) +
O\left(\frac{1}{\l _i}+ \sum \e _{ij}\right),
\end{align*}
where $c_1$,...,$c_6$ are positive constants, and $(e_1,\cdots , e_4) $ denotes an orthonormal basis of $T_{a_i}S^4_+$.
\end{pro}

\begin{pro}\label{p:26}
For $ u=\sum_{i\leq q}\a _i\d _i+\sum_{j=q+1}^p\a _j\var _j \in V(p,\e )$, we 
have the following expansions 
\begin{align*}
(\n J(u),  \l _j\partial\var _j/\partial\l _j)=&  2 \o_3J(u)\left(-2\sum_{k\ne j}\a_k\l_j\frac{\partial \e _{jk}}{\partial \l _j}+4\sum_{k=q+1,k\ne j}^p\a_k\frac{H(a_j,a_k)}{\l _j\l_k}\right.\\
 &  +\a_j\frac{\D K(a_j)}{3\l_j^2K(a_j)}+4\a_j\frac{H(a_j,a_j)}{\l_j^2}
 +o\biggl(\sum_{k=1}^p\frac{1}{(\l_kd_k)^2}+\sum_{k\ne r}\e_{kr}\biggr)\\
(\n J(u),\l _j^{-1}\partial\var _j/\partial a _j).&  \n K(a_j) \geq  c\frac{|\n K(a_j)|^2}{\l_j}+O\biggl(\frac{1}{(\l_jd_j)^2}+\sum_{k\ne j}\e_{kj}\biggr)\\
 3\int_{S_+^4} K\var_i\var_j ^2\l_j\frac{\partial\var_j}{\partial\l_j}& =\frac{\o_3}{4}K(a_j)\left(\frac{\partial\e _{ij}}{\partial\l_j}-2\frac{H(a_i,a_j)}{\l_i\l_j}\right)+o\left(\e_{ij}+\sum\frac{1}{(\l_kd_k)^2}\right)
\end{align*}
\end{pro}
\begin{pf}
First observe that   easy computations show the following estimates:
 
\begin{align*}
(\d_i,\l_j\frac{\partial\var_j}{\partial\l_j}) & =2\o_3\l_j\frac{\partial\e_{ij}}{\partial\l_j}+o(\e_{ij}+\frac{1}{(\l_jd_j)^2})\\
(\var_i,\l_j\frac{\partial\var_j}{\partial\l_j}) & =2 \o_3\l_j\frac{\partial\e_{ij}}{\partial\l_j}-4\o_3\frac{H(a_i,a_j)}{\l_i\l_j}+  o(\e_{ij}+\sum\frac{1}{(\l_jd_j)^2})\\
(\var_j,\l_j\frac{\partial\var_j}{\partial\l_j}) & = -4\o_3\frac{H(a_j,a_j)}{\l_j^2}+o(\frac{1}{(\l_jd_j)^2})\\
\int_{S_+^4}K\d_i ^3\l_j\frac{\partial\var_j}{\partial\l_j} & =\frac{\o_3}{4}K(a_i)\frac{\partial\e_{ij}}{\partial\l_j}+o(\e_{ij}+\frac{1}{(\l_jd_j)^2})\\
\int_{S_+^4}K\var_i ^3\l_j\frac{\partial\var_j}{\partial\l_j} & =\frac{\o_3}{4}K(a_i)\left(\l_j\frac{\partial\e_{ij}}{\partial\l_j}-2\frac{H(a_i,a_i)}{\l_j\l_i}\right)+o(\e_{ij}+\sum\frac{1}{(\l_kd_k)^2})\\
\int_{S_+^4}K\var_j ^3\l_j\frac{\partial\var_j}{\partial\l_j} & =-\frac{\o_3}{24}\frac{\D K(a_i)}{\l_j ^2}- w_3K(a_j)\frac{H(a_j,a_j)}{\l_j ^2}+o(\sum\frac{1}{(\l_jd_j)^2})\\
3\int_{S_+^4}K\d_i \var_j^2\l_j\frac{\partial\var_j}{\partial\l_j} & =\frac{\o_3}{4}K(a_j)\frac{\partial\e_{ij}}{\partial\l_j}+o(\e_{ij})\\
3\int_{S_+^4}K\var_i \var_j^2\l_j\frac{\partial\var_j}{\partial\l_j} & =\frac{\o_3}{4}K(a_j)\left(\l_j\frac{\partial\e_{ij}}{\partial\l_j}-2\frac{H(a_i,a_i)}{\l_j\l_i}\right)+o(\e_{ij}+\sum\frac{1}{(\l_kd_k)^2})
\end{align*}
Using the above estimates and the fact that $J(u)^2\a_i ^2 K(a_i)=8+o(1)$, the Proposition  follows using similar arguments as in \cite{BCCH}.
\end{pf}

\section{Construction of a pseudogradient flow}
\mbox{}

In this section we are going to construct a global pseudogradient flow for the Functional $J$ under assumption $(C)$ on $K$. Along its flow lines there can be only finitely many isolated blow up points. Such a flow is defined by combining two basic facts. On the one hand, the first one comes from the Morse Lemma at infinty which moves  points and concentrations as follows: points move according to $ - \nabla K$ if they are interior points and along $\partial_\nu K$ if they are boundary points, concentrations move so as to decrease the Functional $J$. On the other hand, there is another pseudogradient when the points are very close and the total interaction $\sum \e_{ij}$ is large with respect to $\sum \frac{1}{\l_i^2}$. We need to convex-combine both flows to keep the pseudogradient property, to avoid the creation of new asymptotes and to ensure the property that the flow lines when they leave some $V(p,\e)$ will loose at least one bubble, that is the flow will never come back to $V(q,\e)$ for $q \leq p$, a fact which is not trivial in scalar curvature problems whose functional's levels on $V(p,\e)$ are not constant. Some levels of $V(p,\e)$ might be below some other levels of $V(q,\e)$ for some $ q < p$.

As a by product of the construction of our pseudogradient, we able to   identify   the critical points at infinity of our problem. We recall that the critical point at infinity are  the orbits of the gradient flow of $J$ which remain in $V(p,\e(s))$, where $\e(s)$ , a given function, tends to zero when $s$ tends to $+\infty$ (see \cite{B1}). \begin{pro}\label{p:31}
For $p\geq 1$, there exists a pseudo-gradient $W$ so that the following holds:\\ There is a constant $c > 0$ independent of $u=\sum_{i=1}^q\a_i\d_i+\sum_{j=q+1}^p\a_j\var_j\in V(p,\e)$ so that 
$$
(-\n J(u),W)\geq c\biggl(\sum_{k\ne r}\e_{kr}+\sum_{i\leq q}\frac{1}{\l_i}+\frac{|\n K(a_i)|}{\l_i}+\sum_{j=q+1}^p\frac{1}{(\l_jd_j)^2}\biggr)\leqno{ (i)}
$$
$$
(-\n J(u+\ov{v}),W+\frac{\partial\ov{v}}{\partial (\a_i,a_i,\l_i)}(W))\geq c\biggl(\sum_{k\ne r}\e_{kr}+\sum_{i\leq q}\frac{1}{\l_i}+\frac{|\n K(a_i)|}{\l_i}+\sum_{j=q+1}^p\frac{1}{(\l_jd_j)^2}\biggr)\leqno{ (ii)}
$$ 
(iii) $|W|$ is bounded. Furthermore, the only case where the maximum of the $\l_i$'s is not bounded is when each point $a_i$ is close to a critical point $y_{j_i}$ of $K$ with $j_i\ne j_k$ for $i\ne k$ and $\rho(y_{i_1},...,y_{i_p})>0$, where  $\rho(y_{i_1},...,y_{i_p})$ denotes the least eigenvalue of $M(y_{i_1},...,y_{i_p})$. 
\end{pro}
\begin{pf}
Without loss of generality, we can assume that
$$
\l_1\leq ...\leq \l_q, \qquad \mbox{and }\qquad \l_{q+1}\leq ...\leq \l_p.
$$
Let
$$
u=\sum_{i\leq q}\a_i\d_i + \sum_{i>q}\a_i\var_i \in V(p,\e )
$$
Since $(\partial K(z))/(\partial\nu) < 0$ for any critical point $z$ of $K_1=K_{/\partial S^4_+}$, there exist $\mu >0$ and $c>0$ such that 
$(\partial K(a))/(\partial\nu) < -c<0$ for any $a\in \cup B(z_k,2\mu )$. In the first step, we will build a vector field $Y_b$ using the indices $i\in\{1,...,q\}$. For this purpose, we introduce the following sets
\begin{align*}
P&=\{i\leq q/ \exists \mbox{ a sequence } i_1=i,...,i_j \mbox{ s.t. } a_{i_j} \in\cup B(z_k,\mu ) \mbox{ and } d(a_{i_k},a_{i_{k-1}})<\frac{1}{p}\mu \, \forall k\leq j\}\\
I&=\{i\leq q/\l_i\leq M\l_1\}, \quad \mbox{ where } M \mbox{ is a large constant}.\end{align*}
We observe that, if $i\in P$ then there exists $r$ such that $d(a_i,z_r)<2\mu$ and therefore $(\partial K(a_i))/(\partial\nu) < -c<0$. But if $i\notin P$ then, for any $r$ we have $d(a_i,z_r) \geq \mu$. Furthermore if $j\notin P$ and $i\in P$ then $d(a_i,a_j)\geq \frac{\mu}{p}$. We also have for $\l_i\geq \l_j$
\begin{eqnarray}\label{e:31}
-2\l_i\frac{\partial\e_{ij}}{\partial\l_i}-\l_j\frac{\partial\e_{ij}}{\partial\l_j} = 2\e_{ij}(1-2\frac{\l_j}{\l_i}\e_{ij}) + \e_{ij}(1-2\frac{\l_i}{\l_j}\e_{ij})
\geq \e_{ij}(1+o(1))
\end{eqnarray}
Now, we define the following vector fields
$$
Z_1=\sum_{i\in P}\l_i\frac{\partial\d_i}{\partial\l_i}2^i\quad\mbox{and }\quad  Z_2=\sum_{i\in I\backslash P}\frac{1}{\l_i}\frac{\partial\d_i}{\partial a_i}\frac{\n_{T}K(a_i)}{|\n_{T}K(a_i)|}
$$
Using \eqref{e:31} and Proposition \ref{p:25}, we derive
\begin{eqnarray}\label{e:32}
\biggl(-\n J(u),Z_1\biggr) \geq c\sum_{i\in P, j\leq q}\e_{ij}
+ O\biggl(\sum_{i\in P, j> q}\e_{ij}\biggr) + \sum_{i\in P} \frac{c}{\l_i}
 + o\biggl(\sum_{ k,r\leq q}\e_{kr}\biggr)
 + O\biggl(\sum_{k\in P} \frac{1}{\l_k^2}\biggr)
\end{eqnarray}
and
\begin{align}\label{e:33}
\biggl(-\n J(u),Z_2\biggr) &\geq\sum_{i\in I\backslash P}\biggl(\frac{c}{\l_i} + O\biggl(\sum_{j\in I\cup P}\frac{1}{\l_i}|\frac{\partial\e_{ij}}{\partial a_i}|\biggr)\biggr) + O\biggl(\sum_{j\notin I\cup P}\e_{ij} + \sum_{j>q}\e_{ij}\biggr)\notag\\
& +o\biggl(\sum_{k,r\leq q}\e_{kr}\biggr)
+O\biggl(\sum_{k\leq q}\frac{1}{\l_k^2}\biggr)
\end{align}
Notice that for $i,j\in I$,  $\l_id(a_i,a_j)$ is very large and thus
\begin{eqnarray}\label{e:34}
\frac{1}{\l_i}|\frac{\partial\e_{ij}}{\partial a_i}|\leq c\l_jd(a_i,a_j) \e_{ij}^2 \leq \frac{\e_{ij}}{\l_id(a_i,a_j)}=o(\e_{ij})
\end{eqnarray}
We also notice that if $i\notin P$ and $j\in P$, thus $d(a_i,a_j)>\frac{\mu}{p}$ and therefore
\begin{eqnarray}\label{e:35}
 \frac{1}{\l_i}|\frac{\partial\e_{ij}}{\partial a_i}| = O\biggl(\frac{1}{\l_i^3}+\frac{1}{\l_j^3}\biggr)
\end{eqnarray}
Let us define now the following vector field
$$
Z_3=-\sqrt {M} \sum_{i\notin I\cup P}\l_i\frac{\partial\d_i}{\partial\l_i}2^i - m \sum_{i\in I}\l_i\frac{\partial\d_i}{\partial\l_i}, 
$$
where $m$ is a small positive constant.\\
Using \eqref{e:31} and Proposition \ref{p:25}, we derive
\begin{align}\label{e:36}
\biggl(-\n J(u),Z_3\biggr) \geq \, & c \, \sqrt{M}\sum_{i\notin I\cup P}\left( \sum_{j\leq q, j\notin P}\e_{ij} +O\left(\sum_{j\in P} \e_{ij}+\sum_{j>q}\e_{ij} + \sum_{k\leq q}\frac{1}{\l_k^2}+\frac{1}{\l_i}\right)\right)\notag\\
 & +cm\sum_{i\in I}\biggl(\sum_{j\in I}\e_{ij}+ O\biggl(\sum_{j\notin I}\e_{ij}+\frac{1}{\l_i}\biggr)\biggr)
\end{align}
Observe that, if $i\notin I$ and $i\leq q$, we have $\sqrt{M}/\l_i \leq (\l_1\sqrt{M})^{-1} = o(\l_1^{-1})$, for $M$ large enough.\\
We also define the following vector field
$$
Z_4=\sum_{i\leq q} \psi \biggl(\l_i(1- J(u)^2\a_i^2K(a_i))\biggr)\d_i,
$$
where $\psi$ is a $C^\infty$ function which satisfies
$$
\psi (t)=-1 \quad \mbox{if } t>2, \quad \psi (t)=0 \quad \mbox{if }|t|\leq 1 \quad\mbox{and } \psi (t)=1 \quad\mbox{if } t<-2
$$
Using Proposition \ref{p:25}, we derive
\begin{align}\label{e:37}
<-\n J(u),Z_4> \, \geq \sum_{i\leq q} |\psi (\l_i(1-J(u)^2\a_i^2K(a_i)))|\biggl[|1&-J(u)^2\a_i^2K(a_i)|\notag\\
&+O\biggl(\frac{1}{\l_i}+\sum_{j\ne i}\e_{ij}\biggr)\bigg]
\end{align}
Now, we introduce the following vector field
$$
Y_b=M_1Z_1+\sqrt{M_1}Z_2+Z_3+Z_4
$$
Using \eqref{e:32},..., \eqref{e:37}, we derive
\begin{eqnarray}\label{e:38}
<-\n J(u),Y_b> \geq c \sum_{i\leq q}\biggl(\frac{1}{\l_i}+ |1-J(u)^2\a_i^2K(a_i)|\sum_{j\leq q}\e_{ij}\biggr) +O\biggl(\sum_{i\leq q,j>q}\e_{ij}\biggr)
\end{eqnarray}
Secondly, we need to construct a vector field using the indices $i>q$. We claim that, for $i\leq q$ and $j>q$, we have
\begin{eqnarray}\label{e:39}
-\l_j\frac{\partial\e_{ij}}{\partial\l_i}= \e_{ij}\biggl(1-2\frac{\l_i}{\l_j}\e_{ij}\biggr)=\e_{ij}(1+o(1))
\end{eqnarray}
Indeed, if $\l_j\geq \l_i$, our claim is easy, if $\l_j\leq \l_i$, we obtain $\l_id(a_i,a_j)\geq \l_jd(a_i,a_j)>\l_jd(a_j,\partial S^4_+)> \e^{-1}$, thus $ \frac{\l_i}{\l_j}\e_{ij} \leq \frac{1}{\l_jd(a_i,a_j)}=o(1)$, then our claim follows.\\
Now, we see that there exists $d_0>0$ small enough such that for any $a$ satisfying $d_a=d(a,\partial S^4_+) \leq d_0$, we have
$$
H(a,a) \geq M'|\D K(a)| \quad \mbox{and} \quad H(a,a) \backsim \frac{c}{d_a^2},$$
where $ M'$ is a large constant.\\
We need to introduce the subset of points which are close to  the boundary, for that purpose, let us introduce the following set
$$ 
F=\{i>q/\exists \, \,  i_1=i,...,i_l \mbox{ s.t. } d(a_{i_l},\partial S^4_+)< \frac{d_0}{p} \, \mbox{ and } d(a_{i_k},a_{i_{k-1}})<\frac{d_0}{p}\forall \, \, k\leq l\}
$$
It is easy to see that the following claims hold:\\
- if $i\in F$, $j\notin F$, $j\geq q+1$, we have $d(a_i,a_j)\geq d_0/p$\\
- for any $i\in F$, we have $H(a_i,a_i) \geq M' |\D K(a_i)|$ and $H(a_i,a_i)\backsim c/d_0^2$\\
- for any $j\notin F$ and $j\geq q+1$, we have $d_j>d_0/p$. Furthermore, if $i\leq q$, we have $\e_{ij}=o(\l_i ^{-1})$\\
Thus, \eqref{e:38} becomes
\begin{align}\label{e:310}
<-\n J(u),Y_b> \geq & c \sum_{i\leq q}\biggl(\frac{1}{\l_i}+ |1-J(u)^2\a_i^2K(a_i)|+\sum_{k\leq q}\e_{ki}\biggr)\notag\\
& + c\sum_{i\leq q, j\notin F}\e_{ij} + O\biggl(\sum_{k\leq q,r\in F}\e_{kr}\biggr)
\end{align}
Now, we introduce
$$
Z_5 = -\sum_{i\in F} \l_i \frac{\partial\var _i}{\partial\l_i}2^i
$$
Using Proposition \ref{p:25}, we derive
\begin{eqnarray}\label{e:311}
<-\n J(u),Z_5> \geq c \sum_{k\in F} \e_{kr}+c\sum_{k\in F}\frac{1}{(\l_k d_k)^2} +o(\sum\e_{kr})
\end{eqnarray}
Next, we deal with the points which are far away from the boundary.\\
Let $\bar{\l}=\mbox{max} (\l_1,\mbox{min}\{\l_i, i\in F\})$ ( recall that $\l_1=\mbox{min}\{\l_i, i\leq q\}$). \\
We introduce the following sets
\begin{eqnarray}\label{e:*}
L=\{i>q/ i\notin F \mbox{ and } \l_i \leq \bar{\l }/2 \}, \quad L'=\{i>q/i\notin F \mbox{ and } i\notin L\}
\end{eqnarray}
We denote by $Z_6$ the following vector field 
$$
Z_6 = -\sum_{i\in L'} \l_i \frac{\partial\var _i}{\partial\l_i}2^i
$$
Thus, we have
\begin{align}\label{e:312}
<-\n J(u),Z_6> &\geq c \sum_{k\in L'} \e_{kr}+\sum_{k\in L'}O\biggl(\frac{1}{\l_k ^2}\biggl) +o(\sum\e_{kr})\notag\\
&\geq\sum_{k\in L'} \e_{kr}+O\biggl(\frac{1}{\bar{\l} ^2}\biggl) +o(\sum\e_{kr}) \end{align}
Now, we define
$$
Z_7=M_2Z_5+Y_b+Z_6,
$$
where $M_2$ is a positive constant large enough.\\
Using \eqref{e:310},\eqref{e:311} and \eqref{e:312}, we derive
\begin{align}\label{e:313}
<-\n J(u),Z_7> \geq& c \sum_{k,r\in L} \e_{kr}+c\sum_{i\leq q}\biggl(|1-J(u)^2\a_i^2K(a_i|)+\frac{1}{\l_i}\biggr) \notag\\
&+ c\sum_{k\in F\cup L'} \frac{1}{(\l_id_i)^2} +o(\sum\e_{kr}) 
\end{align} 
Now, we observe that if $L\ne \{1,...,p\}$, then 
\begin{eqnarray}\label{e:314}
\mbox{max}\{\l_i, i\in L\} \leq(1/2) \mbox{max}\{\l_i, i=1,...,p\},
\end{eqnarray}
where $L$ is defined in \eqref{e:*}. Notice that, if
 $i\in L$, we have $d_i\geq d_0$, thus the function $H$ and its gradient are bounded. In this case, we can use the vector field ( denoted here $Z_8$) defined in Lemma 3.3 of \cite{BCCH}. We will apply $Z_8$ only to $u_1= \sum_{i\in L}\a_i\var_i$ forgeting the indices $i\notin L$. Thus, we have
\begin{eqnarray}\label{e:315}
<-\n J(u),Z_8(u_1)> \geq c \biggl(\sum_{k,r\in L} \e_{kr}+\sum_{i\in L}\frac{1}{\l_i^2}\biggr) + O\biggl(\sum_{k\in L, r\notin L}\e_{kr}\biggr) +o\biggl(\sum \frac{1}{\l_k^2}\biggr) 
\end{eqnarray} 
Now, we define
$$
Z_9=Z_8+M_3Z_7,
$$
where $M_3$ is a positive constant large enough.\\
Using \eqref{e:313} and \eqref{e:315}, we derive
\begin{eqnarray}\label{e:316}
<-\n J(u),Z_9> \geq c \sum \e_{kr}+\sum_{i\leq q}|1-J(u)^2\a_i^2K(a_i)|\frac{1}{\l_i} + \sum_{i>q} \frac{1}{(\l_id_i)^2} 
\end{eqnarray} 
To obtain the estimate (i), we need to introduce the following vector field
$$
Z_{10}= \sum_{i>q}\frac{1}{\l_i}\frac{\partial\var_i}{\partial a_i} \frac{\n K(a_i)}{|\n K(a_i)|}\psi \left(\l_i |\n K(a_i)|\right),
$$
where $\psi$ is a $C^\infty$ function satisfying $\psi (t)=1$ si $t\geq 2$ and $\psi (t)=0$ si $t\leq 1$.\\
Our vector field $W$ will be the following
$$
W=M_4Z_9+Z_{10},
$$
where $M_4$ is a positive constant large enough.\\
Thus, using \eqref{e:316} and Proposition \ref{p:25}, the estimate (i) follows.\\
Regarding the estimate (ii), it can be obtained once we have (i), using the estimates of $||\n J(u+\bar{v})||$ and $||\bar{v}||$, arguing as in Appendix B of \cite{BCCH}.\\
Now, we observe that if the set $L$ defined in \eqref{e:*} is equal to $\{1,...,p\}$, thus using \eqref{e:314} and the fact that $Z_7$ only decreases the $\l_i$'s, we derive that the maximum of the $\l_i$'s is a decreasing function in this case. In other case, that is $L=\{1,...,p\}$, Claim (ii) follows from the definition of $Z_8$ (see Lemma 3.3 of \cite{BCCH}). Thus, the proof of our proposition is completed.
\end{pf}
\begin{cor}\label{c:32}
Assume that $J$ has no critical points in $\Sig^+$. Then the only critical points at infinity of $J$ correspond to 
$$
\sum_{j=1}^pK(y_{i_j})^{-1/2}\var _{(y_{i_j},\infty)}, \quad\mbox{ with }\,  p\in \N^* \mbox{ and }\rho(y_{i_1},...,y_{i_p}) > 0.
$$ 
Furthermore, such a critical point at infinity has a Morse index equal to\\
 $(5p-1-\sum_{j=1}^p index(K,y_{i_j}))$, where $index(K,y_{i_j})$ is the Morse index of $K$ at $y_{i_j}$. 
\end{cor}
\begin{pf}
From Proposition \ref{p:31}, we know that the only region where the $\l_i$'s are unbounded is when each point $a_i$ is close to a critical point $y_{j_i}$, with $j_i\ne j_k$ for $i\ne k$ and $\rho(y_{i_1},...,y_{i_p}) > 0$. In this region, arguing as in \cite{B2} and \cite{BCCH}, we can find a change of variable
$$
(a_1,...,a_p,\l_1,...,\l_p)\to (\tilde a_1,...,\tilde a_p,\tilde\l_1,...,\tilde\l_p):=(\tilde a,\tilde\l)
$$  
such that 
$$
J\biggl(\sum_{i=1}^p\a_i\var_i+\ov{v}\biggr)=\psi(\a,\tilde a,\tilde\l):=\frac{8S_4^{1/2}\sum_{i=1}^p\a_i ^2}{(\sum_{i=1}^p\a_i ^4K(a_i))^{1/2}}\biggl(1+(c-\eta)\biggl(\sum_{i=1}^p\frac{1}{K(y_{j_i})}\biggr)^{-1}\, {}^t\L M(\tau_p)\L\biggr)
$$
where $\a=(\a_1,...,\a_p)$, $c$ is a positive constant, $\eta$ is a small positive constant, ${}^t\L=(\tilde\l_1,...,\tilde\l_p)$, $\tau_p=(y_{j_i},...,y_{j_p})$ and $S_4=\int_{\R^4}\tilde\d_{(o,1)}^4$.\\
This yields a split of variables $\tilde a$ and $\tilde\l$, thus it is easy to see that if $\tilde a$ is equal to $(y_{j_1},...,y_{j_p})$, only $\tilde\l$ can move. Since $\rho(y_{i_1},...,y_{i_p}) > 0$, in order to decrease the functional $J$, we have to increase $\tilde\l$, and  we obtain a critical point at infinity only in this case.\\
It remains to compute the Morse index of such a critical point at infinity. In order to compute such a Morse index, we observe that $M(\tau_p)$ is definite positive and the function $\psi$ possesses, with respect to the variables $\a_i$'s, an absolute degenerate maximum with one dimensional nullity space. Then the Morse index of such a  critical point at infinity is equal to $(p-1-\sum_{i=1}^p(4-\mbox{index}(K,y_{j_i})))$. Thus our result follows.
\end{pf}
\section{Proof of Theorems}
\mbox{}

In this section we give the proof of Theorems \ref{t:11}and \ref{t:12}.\\
\begin{pfn}{\bf Theorem \ref{t:11} }
For $\eta > 0$ small enough, we introduce the following neighborhood of $\Sig^+$ $$
V_\eta(\Sig^+)=\{u\in \Sig/\, \,e^{2J(u)} J(u)^3|u^-|_{L^4}^2 < \eta \},
$$
where $u^-=\max(0,-u)$.\\
Recall that, we already built in Proposition \ref{p:31} a vector field $W$ defined in $V(p,\e)$ for $p\geq 1$. Outside $\cup_{p\geq 1}V(p,\e/2)$, we will use $-\n J$ and our global vector field $Z$ will be built using a convex combination of $W$ and $-\n J$. $V_\eta(\Sig^+)$ is invariant under the flow line generated by $Z$ (see \cite{BCCH}). Since $V_\eta(\Sig^+)$ is contractible, we have $\chi(V_\eta(\Sig^+))=1$, where $\chi$ is the Euler-Poincare characteristic. Arguing by contradiction, we assume that $J$ has no critical points in $V_\eta(\Sig^+)$. It follows from Corollary \ref{c:32} that the only critical points at infinity of $J$ in $V_\eta(\Sig^+)$ correspond to 
$$
\sum_{j=1}^pK(y_{i_j})^{-1/2}\var _{(y_{i_j},\infty)}, \quad\mbox{ with }\,  p\in \N^* \mbox{ and }\rho(y_{i_1},...,y_{i_p}) > 0,
$$
and such a critical point at infinity has a Morse index equal to  $(5p-1-\sum_{j=1}^p index(K,y_{i_j}))$. \\
Using the vector field $Z$, we have $V_\eta(\Sig^+)$ retracts by deformation on $\cup W_u(w_\infty)$ (see section 7 and 8 of \cite{BR}), where $W_u(w_\infty)$ is the unstable manifold at infinity of a critical point at infinity $w_\infty$. Then, we have 
$$
1=\chi(V_\eta(\Sig^+))=\sum_{p=1}^l\sum_{\tau_p=(i_1,...,i_p)/\rho(\tau_p)>0}(-1)^{(5p-1-\sum_{j=1}^p index(K,y_{i_j}))},
$$
which is in contradiction with the assumption of our theorem. Thus there exists a critical point of $J$ in $V_\eta(\Sig^+)$. Arguing as in \cite{BCCH}, we prove that this critical point is positive and hence our result follows.
\end{pfn}\\

\begin{pfn}{\bf Theorem \ref{t:12} }
Again, we argue by contradiction. We assume that $J$ has no critical points in $V_\eta (\Sig^+)$. We observe that
$$
H(y,y) = \frac{1}{4} d_y^{-2}, \qquad \qquad \mbox{ for any } \quad y\in S^4_+,
$$
where $d_y=d(y,\partial S^4_+)$.\\
Thus, under the assumption of our theorem, we derive that
$$
(-\D K(y)/(3K(y))) - H(y,y) < 0, \qquad \quad \mbox{for any critical point } y \mbox{ of } K.
$$
Using Corollary \ref{c:32}, we deduce that there is no critical points at infinity of $J$. Let $Z$ be the vector field defined in the proof of Theorem \ref{t:11}. Let $u_0 \in \Sig^+$ and let $\eta (s,u_0)$ be the one parameter group generated by $Z$. It is known that $|\n J|$ is lower bounded outside $V(p,\e /2)$, for any $p\in \N^*$ and for $\e$ small enough, by a fixed constant which depends only on $\e$. Thus, the flow line $\eta (s,u_0)$ cannot remain outside of the set  $V(p,\e /2)$. Furthermore, if the flow line travels from  $V(p,\e /2)$ to the boundary of  $V(p,\e )$, $J(\eta (s,u_0))$ will decrease by a fixed constant which depends on $\e$. Then, this travel cannot be repeated in an infinite time. Thus, there exist $p_0$ and $s_0$ such that the flow line enters into  $V(p_0,\e /2)$ and it does not exit from $V(p_0,\e)$. But in $V(p_0,\e )$, by Proposition \ref{p:31}, we know that the maximum $\l_{max}$ of the $\l_i$'s is bounded by $\l_{max}(s_0)$ and therefore $|\n J|$ is lower bounded. Then when $s$ goes to $+\infty$, $J(u(s))$ goes to $-\infty$ and this yields a contradiction, hence our result follows.
\end{pfn}

\par
\par
\par
\par\par
\par\par
\par\par
\par\par
\bigskip
\noindent
{\bf {Mohamed Ben Ayed}} : D{\'e}partement de Math{\'e}matiques, Facult{\'e} 
des Sciences de Sfax, Route Soukra, Sfax, Tunisia. E-mail: \texttt{
Mohamed.Benayed@fss.rnu.tn}.\\

\noindent
{\bf{Khalil El Mehdi}} : Facult\'e des Sciences et Techniques, Universit\'e de Nouakchott, BP 5026, Nouakchott, Mauritania. E-mail: \texttt{khalil@unic-nkc.mr}\\
and \\
The abdus Salam International Centre for Theoretical Physics, Mathematics 
Section, Strada Costiera, II-34014 Trieste, Italy. E-mail : 
\texttt{elmehdik@ictp.trieste.it}.\\
\par\par

\noindent
{\bf{Mohameden Ould Ahmedou}} :  Rheinische Friedrich-Wilhelms-Universitat 
Bonn,\\ Mathematisches Institut, Beringstrasse 4, D-53115 Bonn, Germany .\\
 E-mail:
\texttt{ahmedou@math.uni-bonn.de}

\begin{thebibliography}{99}

%\bibitem{A}
%T. Aubin,
%\emph{Some nonlinear problem in differential geometry,} 
%Springer-Verlag, New York \textbf{} 1997.

\bibitem{agp} A. Ambrosetti ,  J. Garcia Azorero , A. Peral , \emph{Perturbation of
$\Delta u + u^{\frac{(N+2)}{(N-2)}} = 0$, the Scalar Curvature Problem in
$\mathbb{R}^N$ and related topics}, Journal of  Functional  Analysis, \textbf{165} (1999), 117-149.

\bibitem{ALM} A. Ambrosetti , Y.Y. Li , A. Malchiodi, 
\emph{On the Yamabe problem and the scalar curvature problems under boundary conditions, }  Math. Ann. \textbf{322} (2002), 667--699.

\bibitem{A}

T. Aubin,
Some nonlinear problems in differential geometry, 
Springer-Verlag, New York \textbf{} 1997.

 \bibitem{aujfa} T. Aubin , \emph{Meilleures constantes dans le th\'eor\`eme
d' inclusion de Sobolev et un th\'eor\`eme de Fredholm non lin\'eaire pour
la transformation conforme de la courbure scalaire}, Journal of Functional Analysis,
\textbf{32}, 1979, 148-174.

\bibitem{ab} T. Aubin  and A. Bahri   \emph{M\'ethodes de topologie alg\'ebrique
pour le probl\'eme de la courbure scalaire prescrite}, Journal des Math\'ematiques Pures et
Appliqu\'ees,  \textbf{76} (1997), 525-549.

\bibitem{B1}
 A. Bahri,
 Critical point at infinity in some variational problems,
 Pitman Res. Notes Math, Ser \textbf{ 182}, Longman Sci. Tech. Harlow 1989.
 \bibitem{B2}
A. Bahri,
\emph{An invariant for Yamabe-type flows with applications to scalar curvature
problems in high dimension,} A celebration of J. F. Nash Jr., Duke Math. J. 
\textbf{81} (1996), 323-466.
\bibitem{BC'} 
 A. Bahri and J. M. Coron,
\emph{On a nonlinear elliptic equation involving the critical Sobolev exponent: The effect of topology of the domain, }
 Comm. Pure Appl. Math. \textbf{41} (1988), 255-294.
\bibitem{BaC} 
 A. Bahri and J. M. Coron,
\emph{The scalar curvature problem on the standard three dimensional spheres,}
 J. Funct. Anal. \textbf{95} (1991), 106-172.
\bibitem{BR}
A. Bahri and P. Rabinowitz,
\emph{Periodic orbits of hamiltonian systems of three body type, }
Ann. Inst. H. Poincar{\'e} Anal. Non lin{\'e}aire \textbf{8} (1991), 561-649.

\bibitem{BCCH}
 M. Ben Ayed, Y. Chen, H. Chtioui and M. Hammami,
\emph{ On the prescribed scalar curvature problem on 4-manifolds,}
 Duke Math. J. \textbf{84} (1996), 633-677.

\bibitem{BEO}
 M. Ben Ayed, K. El Mehdi  and M. Ould Ahmedou,
\emph{Prescribing the scalar curvature under minimal boundary conditions on the half sphere, } Adv. Nonlinear Stud. \textbf{2} (2002), 93-116.

\bibitem{bp}
 G. Bianchi and X. B. Pan,
 \emph{Yamabe equations on
half-spaces}, Nonlinear Anal. \textbf{37} (1999) 161-186.

\bibitem{BC}
H. Brezis and J. M. Coron,
\emph{Convergence of solutions of H-systems or how to blow bubbles,}
 Arch. Rational Mech. Anal. \textbf{89} (1985), 21-56.
\bibitem{cgy} 
S-Y. A. Chang , M.J. Gursky and P.C. Yang , 
\emph{The scalar
curvature equation on 2- and 3- spheres,}
Calculus of Variations and Partial Differential Equations, \textbf{1} (1993), 205-229.

\bibitem{cy1}
 S-Y.A Chang and P. Yang , \emph{A perturbation result in
prescribing scalar curvature on $S^n$}, Duke Mathematical  Journal, \textbf{64} (1991),27-69.

\bibitem{cl} 
K. C.  Chang and J. Q. Liu ,
\emph{On Nirenberg's problem,} Int. J. Math. 4 (1993), 35-58.
\bibitem{C}
P. Cherrier,
\emph{Probl{\`e}mes de Neumann non lin{\'e}aires sur les vari{\'e}t{\'e}s 
Riemaniennes,}
J. Funct. Anal. \textbf{57} (1984), 154-207.
\bibitem{DMO}
Z. Djadli, A. Malchiodi and M. Ould Ahmedou,
\emph{Prescribing the scalar and the boundary mean curvature on the three 
dimensional half sphere,}  J. Geom. Anal. \textbf{13} (2003), 233 - 267. 
\bibitem{E1}
J. Escobar,
\emph{Conformal deformation of Riemannian metric to scalar flat metric
 with constant  mean curvature on the boundary,}
Ann. of Math. \textbf{136} (1992), 1-50.
\bibitem{E2}
J. Escobar,
\emph{Conformal metrics with prescribed mean curvature on the boundary,}
Cal. Var. \textbf{4} (1996), 559-592.
\bibitem{es}
J.  Escobar , R. Schoen, 
\emph{Conformal metrics
with prescribed scalar curvature}, Inventiones Mathematicae, 86 (1986), 243-254.
\bibitem{Escobar} 
J. Escobar, 
\emph{On the precribed scalar curvature on compact manifolds with boundary}, Differential geometric methods in the control of partial differential equations(Boulder, CO, 1999), Contemp. Math. \textbf{268}(2000), 137 - 144.

\bibitem{HL1}
Z. C. Han and Y.Y. Li,
\emph{The Yamabe problem on manifolds with boundaries : existence and 
compactness results,}
 Duke Math. J. \textbf{99} (1999), 489-542.
\bibitem{HL2}
Z. C. Han and Y.Y. Li,
\emph{The existence of conformal metrics with constant scalar curvature 
and constant boundary mean curvature,}
Comm. Anal. Geom. \textbf{8} (2000), 809-869.

\bibitem{hv} E. Hebey , \emph{Changements de m\'etriques conformes sur la sph\`ere -
Le probl\`eme de Nirenberg}, Bull. Sci. Math. \textbf{114} (1990), 215-242.

\bibitem{Hebey}
E. Hebey, \emph{The isometry concentration method in the case of a nonlinear problem with Sobolev critical exponent on compact manifolds with boundary}, Bull. Sci. Math. \textbf{116}(1992), 35 - 51.
\bibitem{H}
E. Hebey,
Introduction \`a l'analyse non lin\'eaire sur les variet\'es, Diderot Editeur, Paris, 1997.

\bibitem{y1} Li Y.Y., \emph{Prescribing scalar curvature on $S^{n}$ and
related topics, Part I}, Journal of  Differential  Equations,   \textbf{120} (1995), 319-410; \emph{Part
II, Existence and compactness}, Communications in   Pure and Applied  Mathematics,  \textbf{49} (1996),
437-477.



\bibitem{Li}
Y.Y. Li,
\emph{The Nirenberg problem in a domain with boundary,}
Top. Meth. Nonlin. Anal. \textbf{6} (1995), 309-329.
\bibitem{L}
P. L. Lions,
\emph{The concentration compactness principle in the calculus of variations.
 The limit case,}
Rev. Mat. Iberoamericana \textbf{1} (1985), I: 165-201; II: 45-121.
\bibitem{R}
O. Rey,
\emph{Boundary effect for an elliptic Neumann problem with critical 
nonlinearity,}
Comm. Partial Diff. Eq. \textbf{22} (1997), 1055-1139.

\bibitem{[SZ]} R. Schoen , D. Zhang , \emph{Prescribed scalar
curvature on the $n$-sphere}, Calculus of  Variations and Partial Differential Equations,
\textbf{4} (1996), 1-25.

\bibitem{S}
M. Struwe,
\emph{A global compactness result for elliptic boundary value problems 
involving nonlinearities,}
Math. Z. \textbf{187} (1984), 511-517. 
\end{thebibliography}
\end{document}